\documentclass{elsart}
\usepackage{amsmath}
\usepackage{amssymb}

\newcommand{\N}{\mathbb{N}}                     
\newcommand{\R}{\mathbb{R}}                     
\newcommand{\dd}{\, {\rm d}}
\newcommand{\E}{\mathbf{E}}
\newcommand{\med}{{\rm median}}

\newcommand{\sug}{\mathcal{S}_{\mu}}
\newcommand{\cho}{\mathcal{C}_{\mu}}

\begin{document}

\begin{frontmatter}

\title{Weighted lattice polynomials of independent random variables}

\author{Jean-Luc Marichal}
\ead{jean-luc.marichal[at]uni.lu}

\address{
Applied Mathematics Unit, University of Luxembourg\\
162A, avenue de la Fa\"{\i}encerie, L-1511 Luxembourg, Luxembourg }

\date{June 5, 2007}

\begin{abstract}
We give the cumulative distribution functions, the expected values, and the moments of weighted lattice polynomials when regarded as real
functions of independent random variables. Since weighted lattice polynomial functions include ordinary lattice polynomial functions and,
particularly, order statistics, our results encompass the corresponding formulas for these particular functions. We also provide an application
to the reliability analysis of coherent systems.
\end{abstract}

\end{frontmatter}

\section{Introduction}

The cumulative distribution functions (c.d.f.'s) and the moments of order statistics have been discovered and studied for many years (see for
instance David and Nagaraja \cite{DavNag03}). Their generalizations to lattice polynomial functions, which are nonsymmetric extensions of order
statistics, were investigated very recently by Marichal \cite{Marb} for independent variables and then by Dukhovny \cite{Duk} for dependent
variables.

Roughly speaking, an $n$-ary {\em lattice polynomial}\/ is a well-formed expression involving $n$ real variables $x_1,\ldots,x_n$ linked by the
lattice operations $\wedge=\min$ and $\vee=\max$ in an arbitrary combination of parentheses. In turn, such an expression naturally defines an
$n$-ary {\em lattice polynomial function}. For instance,
$$
p(x_1,x_2,x_3)=(x_1\wedge x_2)\vee x_3
$$
is a 3-ary (ternary) lattice polynomial function.

Lattice polynomial functions can be straightforwardly generalized by fixing certain variables as parameters, like in the 2-ary (binary)
polynomial function
$$p(x_1,x_2)=(c\wedge x_1)\vee x_2,$$ where $c$ is a real constant. Such ``parameterized'' lattice polynomial functions, called {\em weighted
lattice polynomial}\/ functions \cite{Marc}, play an important role in the areas of nonlinear aggregation and integration. Indeed, they include
the whole class of discrete Sugeno integrals \cite{Sug74,Sug77}, which are very useful aggregation functions in many areas. More details about
the Sugeno integrals and their applications can be found in the remarkable edited book \cite{GraMurSug00}.

In this paper we give closed-form formulas for the c.d.f.\ of any weighted lattice polynomial function in terms of the c.d.f.'s of its input
variables. More precisely, considering a weighted lattice polynomial function $p:\R^n\to\R$ and $n$ independent random variables
$X_1,\ldots,X_n$, $X_i$ $(i=1,\ldots,n)$ having c.d.f.\ $F_i(x)$, we give formulas for the c.d.f.\ of $Y_p:=p(X_1,\ldots,X_n)$. We also yield
formulas for the expected value $\E[g(Y_p)]$, where $g$ is any measurable function. The special cases $g(x)=x$, $x^r$, $[x-\E(Y_p)]^r$, and
$e^{tx}$ give, respectively, the expected value, the raw moments, the central moments, and the moment-generating function of $Y_p$.

As the lattice polynomial functions are particular weighted lattice polynomial functions, we retrieve, as a corollary, the formulas of the
c.d.f.'s and the moments of the lattice polynomial functions.

This paper is organized as follows. In Section 2 we recall the basic material related to lattice polynomial functions and their weighted
versions. In Section 3 we provide the announced results. In Section 4 we investigate the particular case where the input random variables are
uniformly distributed over the unit interval. Finally, in Section 5 we provide an application of our results to the reliability analysis of
coherent systems.

\section{Weighted lattice polynomials}

In this section we recall some basic definitions and properties related to weighted lattice polynomial functions. More details and proofs can be
found in \cite{Marc}.

As we are concerned with weighted lattice polynomial functions of random variables, we do not consider weighted lattice polynomial functions on
a general lattice, but simply on a closed interval $L:=[a,b]$ of the extended real number system $\overline{\R}:=\R\cup\{-\infty,+\infty\}$.
Clearly, such an interval is a bounded distributive lattice, with $a$ and $b$ as bottom and top elements. The lattice operations $\wedge$ and
$\vee$ then represent the minimum and maximum operations, respectively. To simplify the notation, we also set $[n]:=\{1,\ldots,n\}$ for any
integer $n\geqslant 1$.

Let us first recall the definition of lattice polynomial functions (with real variables); see for instance Birkhoff \cite[\S II.5]{Bir67}.

\begin{defn}
The class of \emph{lattice polynomial functions}\/ from $L^n$ to $L$ is defined as follows:
\begin{enumerate}
\item For any $k\in [n]$, the projection $(x_1,\ldots,x_n)\mapsto x_k$ is a lattice polynomial function from $L^n$ to $L$.

\item If $p$ and $q$ are lattice polynomial functions from $L^n$ to $L$, then $p\wedge q$ and $p\vee q$ are lattice polynomial functions from
$L^n$ to $L$.

\item Every lattice polynomial function from $L^n$ to $L$ is constructed by finitely many applications of the rules (1) and (2).
\end{enumerate}
\end{defn}

As mentioned in the introduction, weighted lattice polynomial functions generalize lattice polynomial functions by considering both variables
and constants. We thus have the following definition.

\begin{defn}\label{de:wlp}
The class of \emph{weighted lattice polynomial functions}\/ from $L^n$ to $L$ is defined as follows:
\begin{enumerate}
\item For any $k\in [n]$ and any $c\in L$, the projection $(x_1,\ldots,x_n)\mapsto x_k$ and the constant function $(x_1,\ldots,x_n)\mapsto c$
are weighted lattice polynomial functions from $L^n$ to $L$.

\item If $p$ and $q$ are weighted lattice polynomial functions from $L^n$ to $L$, then $p\wedge q$ and $p\vee q$ are weighted lattice polynomial
functions from $L^n$ to $L$.

\item Every weighted lattice polynomial function from $L^n$ to $L$ is constructed by finitely many applications of the rules (1) and (2).
\end{enumerate}
\end{defn}

Because $L$ is a distributive lattice, any weighted lattice polynomial function can be written in {\em disjunctive}\/ and {\em conjunctive}\/
forms as follows.

\begin{prop}\label{prop:wlp dnf}
Let $p:L^n\to L$ be any weighted lattice polynomial function. Then there are set functions $\alpha:2^{[n]}\to L$ and $\beta:2^{[n]}\to L$ such
that
$$
p(\mathbf{x})=\bigvee_{S\subseteq [n]}\Big[\alpha(S)\wedge\bigwedge_{i\in S}x_i\Big]=\bigwedge_{S\subseteq [n]}\Big[\beta(S)\vee\bigvee_{i\in
S}x_i\Big].
$$
\end{prop}

Proposition~\ref{prop:wlp dnf} naturally includes the classical lattice polynomial functions. To see this, it suffices to consider nonconstant
set functions $\alpha:2^{[n]}\to\{a,b\}$ and $\beta:2^{[n]}\to\{a,b\}$, with $\alpha(\varnothing)=a$ and $\beta(\varnothing)=b$.

The set functions $\alpha$ and $\beta$ that disjunctively and conjunctively define the polynomial function $p$ in Proposition~\ref{prop:wlp dnf}
are not unique. For example, we have $$x_1 \vee (x_1 \wedge x_2) = x_1=x_1\wedge (x_1 \vee x_2).$$ However, it can be shown that, from among all
the possible set functions that disjunctively define a given weighted lattice polynomial function, only one is nondecreasing. Similarly, from
among all the possible set functions that conjunctively define a given weighted lattice polynomial function, only one is nonincreasing. These
particular set functions are given by
$$
\alpha(S) = p(\mathbf{e}_S)\quad\mbox{and}\quad \beta(S) = p(\mathbf{e}_{[n]\setminus S}),
$$
where, for any $S\subseteq [n]$, $\mathbf{e}_S$ denotes the characteristic vector of $S$ in $\{a,b\}^n$, i.e., the $n$-dimensional vector whose
$i$th component is $b$, if $i\in S$, and $a$, otherwise. Thus, any $n$-ary weighted lattice polynomial function can always be written as
$$
p(\mathbf{x})=\bigvee_{S\subseteq [n]}\Big[p(\mathbf{e}_S)\wedge\bigwedge_{i\in S}x_i\Big]=\bigwedge_{S\subseteq
[n]}\Big[p(\mathbf{e}_{[n]\setminus S})\vee\bigvee_{i\in S}x_i\Big].
$$

The best known instances of weighted lattice polynomial functions are given by the discrete {\em Sugeno integrals}, which consist of a nonlinear
discrete integration with respect to a {\em fuzzy measure}.

\begin{defn}
An $L$-valued {\em fuzzy measure}\/ on $[n]$ is a nondecreasing set function $\mu:2^{[n]}\to L$ such that $\mu(\varnothing)=a$ and $\mu([n])=b$.
\end{defn}

The Sugeno integrals can be presented in various equivalent forms. The next definition introduce them in one of their simplest forms (see Sugeno
\cite{Sug74}).

\begin{defn}
Let $\mu$ be an $L$-valued fuzzy measure on $[n]$. The {\em Sugeno integral}\/ of a function $\mathbf{x}:[n]\to L$ with respect to $\mu$ is
defined by
$$
\sug(\mathbf{x}):=\bigvee_{S\subseteq [n]}\Big[\mu(S)\wedge\bigwedge_{i\in S}x_i\Big].
$$
\end{defn}

Thus, any function $f:L^n\to L$ is an $n$-ary Sugeno integral if and only if it is a weighted lattice polynomial function fulfilling
$f(\mathbf{e}_{\varnothing})=a$ and $f(\mathbf{e}_{[n]})=b$.

\section{Cumulative distribution functions and moments}

Consider $n$ independent random variables $X_1,\ldots,X_n$, $X_i$ $(i\in [n])$ having c.d.f.\ $F_i(x)$, and set $Y_p:=p(X_1,\ldots,X_n)$, where
$p:L^n\to L$ is any weighted lattice polynomial function. Let $H:\overline{\R}\to\{0,1\}$ be the Heaviside step function, defined by $H(x)=1$,
if $x\geqslant 0$, and $0$, otherwise. For any $c\in \overline{\R}$, we also define the function $H_c(x):=H(x-c)$.

The c.d.f.\ of $Y_p$ is given in the next theorem.

\begin{thm}\label{thm:main}
Let $p:L^n\to L$ be a weighted lattice polynomial function. Then, the c.d.f.\ of $Y_p$ is given by each of the following formulas:
\begin{eqnarray}
F_p(y)&=& 1-\sum_{S\subseteq [n]} \big[1-H_{p(\mathbf{e}_S)}(y)\big]\prod_{i\in [n]\setminus
S}F_i(y)\,\prod_{i\in S}[1-F_i(y)],\label{eq:fpy1}\\
F_p(y)&=& \sum_{S\subseteq [n]}H_{p(\mathbf{e}_{[n]\setminus S})}(y)\,\prod_{i\in S}F_i(y)\,\prod_{i\in [n]\setminus
S}[1-F_i(y)].\label{eq:fpy2}
\end{eqnarray}
\end{thm}

\begin{pf*}{Proof.}
Starting from the disjunctive form of $p$, we have
\begin{eqnarray*}
F_p(y) &=& 1-\Pr\Big[\bigvee_{S\subseteq [n]}\big[p(\mathbf{e}_S)\wedge\bigwedge_{i\in S}X_i\big] > y\Big]\\
&=& 1-\Pr\Big[\exists S\subseteq [n]\mbox{ such that }y < p(\mathbf{e}_S)~\mbox{and}~y < \bigwedge_{i\in S}X_i\Big].
\end{eqnarray*}

Consider the following events:
\begin{eqnarray*}
A &:=& \Big[\exists S\subseteq [n]~\mbox{such that}~y < p(\mathbf{e}_S)~\mbox{and}~y < \bigwedge_{i\in S}X_i\Big],\\
B &:=& \Big[\exists S\subseteq [n]~\mbox{such that}~y < p(\mathbf{e}_S)~\mbox{and}~\bigvee_{i\in [n]\setminus S}X_i\leqslant y < \bigwedge_{i\in
S}X_i\Big].
\end{eqnarray*}
Event $B$ implies event $A$ trivially. Conversely, noting that $p$ is nondecreasing in each variable and replacing $S$ with a superset
$S'\supseteq S$, if necessary, we readily see that event $A$ implies event $B$.

Since the events $\big[\bigvee_{i\in [n]\setminus S}X_i\leqslant y < \bigwedge_{i\in S}X_i\big]$ ($S\subseteq [n]$) are mutually exclusive, we
have
$$
F_p(y) = 1- \sum_{S\subseteq [n]}\Pr\big[y < p(\mathbf{e}_S)\big]\,\Pr\Big[\bigvee_{i\in [n]\setminus S}X_i\leqslant y < \bigwedge_{i\in
S}X_i\Big],
$$
which, using independence, proves the first formula. The second one can be proved similarly by starting from the conjunctive form of $p$.\qed
\end{pf*}

The expressions of $F_p(y)$, given in Theorem~\ref{thm:main}, are closely related to the following concept of {\em multilinear extension}\/ of a
set function, which was introduced by Owen \cite{Owe72} in game theory.

\begin{defn}
The {\em multilinear extension}\/ of a set function $v:2^{[n]}\to\R$ is the function $\Phi_v:[0,1]^n\to\R$ defined by
$$
\Phi_v(\mathbf{x}):=\sum_{S\subseteq [n]}v(S)\prod_{i\in S}x_i\, \prod_{i\in [n]\setminus S}(1-x_i).
$$
\end{defn}

Using this concept, we can immediately rewrite (\ref{eq:fpy1}) and (\ref{eq:fpy2}) as
\begin{eqnarray*}
F_p(y)&=& 1-\Phi_{v_{p,y}}[1-F_1(y),\ldots,1-F_n(y)],\\
F_p(y)&=& \Phi_{v^*_{p,y}}[F_1(y),\ldots,F_n(y)],
\end{eqnarray*}
where, for any fixed $y\in\R$, the (nondecreasing) set functions $v_{p,y}:2^{[n]}\to\{0,1\}$ and $v^*_{p,y}:2^{[n]}\to\{0,1\}$ are defined by
$$
v_{p,y}(S):=1-H_{p(\mathbf{e}_S)}(y) \quad\mbox{and}\quad v^*_{p,y}(S):=H_{p(\mathbf{e}_{[n]\setminus S})}(y).
$$

Owen \cite{Owe72} showed that the function $\Phi_v$, being a multilinear polynomial, has the form
\begin{equation}\label{eq:MultPolMoe}
\Phi_v(\mathbf{x})=\sum_{S\subseteq [n]} m_v(S)\,\prod_{i\in S}x_i,
\end{equation}
where the set function $m_v:2^{[n]}\to\R$, called the {\em M\"obius transform}\/ of $v$, is defined as
\begin{equation}\label{eq:motr}
m_v(S)=\sum_{T\subseteq S}(-1)^{|S|-|T|} v(T).
\end{equation}

Using this polynomial form of $\Phi_v$, we can immediately derive two further formulas for $F_p(y)$, namely
\begin{eqnarray}
F_p(y)&=& 1-\sum_{S\subseteq [n]} m_{v_{p,y}}(S)\, \prod_{i\in S}[1-F_i(y)],\label{eq:fpy3}\\
F_p(y)&=& \sum_{S\subseteq [n]} m_{v^*_{p,y}}(S)\, \prod_{i\in S}F_i(y).\label{eq:fpy4}
\end{eqnarray}

Formulas (\ref{eq:fpy1})--(\ref{eq:fpy2}) and (\ref{eq:fpy3})--(\ref{eq:fpy4}) thus provide four equivalent expressions for $F_p(y)$. As
particular cases, we retrieve the c.d.f.\ of any lattice polynomial function. For example, using formula (\ref{eq:fpy1}) leads to the following
corollary (see \cite{Marb}).

\begin{cor}\label{cor:fosk1}
Let $p:L^n\to L$ be a lattice polynomial function. Then, the c.d.f.\ of $Y_p$ is given by
$$
F_p(y)=1-\sum_{\textstyle{S\subseteq [n]\atop p(\mathbf{e}_S)=b}} \prod_{i\in [n]\setminus S}F_i(y)\,\prod_{i\in S}[1-F_i(y)].
$$
\end{cor}

Let us now consider the expected value $\E[g(Y_p)]$, where $g:\overline{\R}\to\overline{\R}$ is any measurable function. For instance, when
$g(x)=x^r$, we obtain the raw moments of $Y_p$.

By definition, we simply have
$$
\E[g(Y_p)]=\int_{-\infty}^{\infty}g(y)\dd F_p(y)=-\int_{-\infty}^{\infty}g(y)\dd [1-F_p(y)].
$$

Using integration by parts, we can derive alternative expressions of $\E[g(Y_p)]$. We then have the following immediate result.

\begin{prop}\label{prop:egydg}
Let $p:L^n\to L$ be any weighted lattice polynomial function and let $g:\overline{\R}\to\overline{\R}$ be any measurable function of bounded
variation.
\begin{enumerate}
\item If $\lim_{y\to\infty}g(y)[1-F_i(y)]=0$ for all $i\in [n]$, then
$$
\E[g(Y_p)]= g(-\infty)+\int_{-\infty}^{\infty}[1-F_p(y)] \dd g(y).
$$
\item If $\lim_{y\to -\infty}g(y)F_i(y)=0$ for all $i\in [n]$, then
$$
\E[g(Y_p)]=g(\infty)-\int_{-\infty}^{\infty}F_p(y) \dd g(y).
$$
\end{enumerate}
\end{prop}

Clearly, combining Proposition~\ref{prop:egydg} with formulas (\ref{eq:fpy1})--(\ref{eq:fpy2}) and (\ref{eq:fpy3})--(\ref{eq:fpy4}) immediately
leads to various explicit expressions of $\E[g(Y_p)]$. For instance, if $$\lim_{y\to\infty}g(y)[1-F_i(y)]=0, \qquad \forall\, i\in [n],$$ then
\begin{equation}\label{eq:egyp}
\E[g(Y_p)]=g(-\infty)+\sum_{S\subseteq [n]}\int_{-\infty}^{p(\mathbf{e}_S)}\prod_{i\in [n]\setminus S}F_i(y)\,\prod_{i\in S}[1-F_i(y)] \dd g(y).
\end{equation}

It is noteworthy that Eq.\ (\ref{eq:egyp}) can also be established without the knowledge of the c.d.f.\ of $Y_p$. As the proof is very
informative, we give it in the appendix.

\begin{rem}
We can also retrieve the c.d.f.\ of $Y_p$ directly from formula (\ref{eq:egyp}). Indeed, rewriting (\ref{eq:egyp}) with the function
$g(y)=H(z-y)$, we immediately obtain
\begin{eqnarray*}
F_p(z) &=& \Pr[z-Y_p\geqslant 0] ~ = ~ \Pr[H(z-Y_p)=1] ~ = ~ \E[H(z-Y_p)]\\
&=& 1+\int_{-\infty}^{\infty}\sum_{S\subseteq [n]} \big[1-H_{p(\mathbf{e}_S)}(y)\big]\prod_{i\in [n]\setminus S}F_i(y)\,\prod_{i\in S}[1-F_i(y)]
\dd H(z-y)
\end{eqnarray*}
and hence we retrieve (\ref{eq:fpy1}).
\end{rem}


\begin{rem}
Suppose that some variables $X_i$ are constants, say, $X_k=c_k$ for all $k\in K$ for a given $K\subseteq [n]$. Then, the weighted lattice
polynomial function $p$ reduces to an $(n-|K|)$-ary weighted lattice polynomial function and it is easy to see that
$$
\bigvee_{S\subseteq [n]}\Big[p(\mathbf{e}_S)\wedge\bigwedge_{i\in S}X_i\Big]=\bigvee_{S\subseteq [n]\setminus
K}\Big[\alpha^K_p(S)\wedge\bigwedge_{i\in S}X_i\Big]
$$
where $$\alpha^K_p(S):=\bigvee_{T\subseteq K}\Big[p(\mathbf{e}_{S\cup T})\wedge\bigwedge_{j\in T}c_j\Big]\qquad (S\subseteq [n]\setminus K).$$
Thus, the conditional expectation $\E[g(Y_p)\mid X_k=c_k\,\forall k\in K]$ can be immediately calculated by Proposition~\ref{prop:egydg}.
\end{rem}

\section{The case of uniformly distributed variables on the unit interval}

We now examine the case where the random variables $X_1,\ldots,X_n$ are uniformly distributed on $[0,1]$. We also assume $L=[0,1]$.

Recall that the {\em incomplete beta function}\/ is defined, for any $u,v>0$, by
$$
B_z(u,v):=\int_0^z t^{u-1}(1-t)^{v-1}\dd t\qquad (z\in\R),
$$
and the {\em beta function}\/ is defined, for any $u,v>0$, by $B(u,v):=B_1(u,v)$.

According to Eq.\ (\ref{eq:egyp}), for any weighted lattice polynomial function $p:[0,1]^n\to [0,1]$ and any measurable function
$g:[0,1]\to\overline{\R}$ of bounded variation, we have
$$
\E[g(Y_p)]=g(0)+\sum_{S\subseteq [n]}\int_{0}^{p(\mathbf{e}_S)}y^{n-|S|}(1-y)^{|S|} \dd g(y).
$$
If, furthermore, $g(y)=\frac 1r y^r$, with $r\in\N\setminus\{0\}$, then
\begin{equation}\label{eq:reypr}
\frac 1r\,\E[Y_p^r]=\sum_{S\subseteq [n]}B_{p(\mathbf{e}_S)}(n-|S|+r,|S|+1).
\end{equation}

\begin{rem}
Considering Eq.~(\ref{eq:reypr}), with $p(\mathbf{e}_S)=z$, if $S=[n]$, and $0$, otherwise, we obtain
$$
B_z(r,n+1)=\frac 1r\int_{[0,1]^n}\Big[z\wedge\bigwedge_{i\in [n]}x_i\Big]^r\dd\mathbf{x}
$$
and hence the following identity
$$
\frac 1r\int_{[0,1]^n}\Big(\bigvee_{S\subseteq [n]}\Big[p(\mathbf{e}_S)\wedge\bigwedge_{i\in S}x_i\Big]\Big)^r\dd\mathbf{x} =\sum_{S\subseteq
[n]}\frac 1{n-|S|+r}\int_{[0,1]^n}\Big[p(\mathbf{e}_S)\wedge\bigwedge_{i\in S}x_i\Big]^{n-|S|+r}\dd\mathbf{x},
$$
which shows that computing the raw moments of any weighted lattice polynomial reduces to computing the raw moments of bounded minima.
\end{rem}

Let us now examine the case of the Sugeno integrals. As these integrals are usually considered over the domain $[0,1]^n$, we naturally calculate
their expected values when their input variables are uniformly distributed over $[0,1]^n$. Since any Sugeno integral is a particular weighted
lattice polynomial, by Eq.~(\ref{eq:reypr}), its expected value then writes
\begin{eqnarray*}
\int_{[0,1]^n}\sug(\mathbf{x})\dd\mathbf{x} &=& \sum_{S\subseteq [n]}B_{\mu(S)}(n-|S|+1,|S|+1)\\
&=& \sum_{S\subseteq [n]} \int_0^{\mu(S)}x^{n-|S|}(1-x)^{|S|}\dd\mathbf{x}\\
&=& \sum_{S\subseteq [n]} \sum_{i=0}^{|S|}{|S|\choose i}(-1)^i\,\frac{\mu(S)^{n-|S|+i+1}}{n-|S|+i+1}.
\end{eqnarray*}

Surprisingly, this expression is very close to that of the expected value of the Choquet integral with respect to the same fuzzy measure.

Let us recall the definition of the Choquet integrals \cite{Cho53}. Just as for the Sugeno integrals, the Choquet integrals can be expressed in
various equivalent forms. We present them in one of their simplest forms (see for instance \cite{Mar02b}).

\begin{defn}
Let $\mu$ be a $[0,1]$-valued fuzzy measure on $[n]$. The {\em Choquet integral}\/ of a function $\mathbf{x}:[n]\to [0,1]$ with respect to $\mu$
is defined by
$$
\cho(\mathbf{x}):=\sum_{S\subseteq [n]}m_{\mu}(S)\bigwedge_{i\in S}x_i,
$$
where $m_{\mu}:2^{[n]}\to\R$ is the M\"obius transform (cf.\ (\ref{eq:motr})) of $\mu$.
\end{defn}

For comparison purposes, we recall the expected value of $\cho$ (see for instance \cite{Mar04}):
\begin{eqnarray*}
\int_{[0,1]^n}\cho(\mathbf{x})\dd\mathbf{x} &=& \sum_{S\subseteq [n]} \mu(S)\, B(n-|S|+1,|S|+1)\\
&=& \sum_{S\subseteq [n]} \mu(S)\,\int_0^{1}x^{n-|S|}(1-x)^{|S|}\dd\mathbf{x}\\
&=& \sum_{S\subseteq [n]} \mu(S)\,\frac{(n-|S|)!\, |S|!}{(n+1)!}.
\end{eqnarray*}

\section{Application to reliability theory}

In this final section we show how the results derived here can be applied to the reliability analysis of coherent systems. For a reference on
reliability theory, see for instance Barlow and Proschan \cite{BarPro81}.

Consider a system made up of $n$ independent components, each component $C_i$ ($i\in [n]$) having a lifetime $X_i$ and a reliability
$r_i(t):=\Pr[X_i>t]$ at time $t>0$. Additional components, with constant lifetimes, may also be considered.

We assume that, when components are connected in series, the lifetime of the subsystem they form is simply given by the minimum of the component
lifetimes. Similarly, for a parallel connection, the subsystem lifetime is the maximum of the component lifetimes.

It follows immediately that, for a system mixing series and parallel connections, the system lifetime is given by a weighted lattice polynomial
function $$Y_p=p(X_1,\ldots,X_n)$$ of the component lifetimes. Our results then provide explicit formulas for the c.d.f., the expected value,
and the moments of the system lifetime.

For example, the system reliability at time $t>0$ is given by
\begin{equation}\label{eq:SystRel}
R_p(t):=\Pr[Y_p>t]=\Phi_{v_{p,t}}[r_1(t),\ldots,r_n(t)].
\end{equation}
Moreover, for any measurable function $g:[0,\infty]\,\to\overline{\mathbb{R}}$ of bounded variation and such that
$\lim_{y\to\infty}g(y)r_i(y)=0$ for all $i\in [n]$, we have, by Proposition~\ref{prop:egydg},
$$
\E[g(Y_p)]=g(0)+\int_0^{\infty}R_p(t)\, {\rm d} g(t).
$$

As an example, the following proposition yields the {\em mean time to failure}\/ $\E[Y_p]$ in the special case of the exponential reliability
model.

\begin{prop}\label{prop:exprel}
If $r_i(t)=e^{-\lambda_i t}$ for all $i\in [n]$, we have
$$
\E[Y_p]=p(\mathbf{e}_{\varnothing})+\sum_{\textstyle{S\subseteq [n]\atop S\neq\varnothing}} \, \sum_{T\subseteq S} (-1)^{|S|-|T|}\,
\frac{1-e^{-\lambda(S)\, p(\mathbf{e}_T)}}{\lambda(S)}\, ,
$$
where $\lambda(S):=\sum_{i\in S}\lambda_i$.
\end{prop}

\begin{pf*}{Proof.}
Using (\ref{eq:SystRel}) and then (\ref{eq:MultPolMoe})--(\ref{eq:motr}), we obtain
$$
R_p(t) = \sum_{S\subseteq [n]} m_{v_{p,t}}(S)\, e^{-\lambda(S)t} = \sum_{S\subseteq [n]}\sum_{T\subseteq S} (-1)^{|S|-|T|}\, v_{p,t}(T)\,
e^{-\lambda(S)t}
$$
and hence
$$
\E[Y_p]=\int_0^{\infty}R_p(t)\, {\rm d}t=\sum_{S\subseteq [n]} \, \sum_{T\subseteq S} (-1)^{|S|-|T|}\,
\int_0^{p(\mathbf{e}_T)}e^{-\lambda(S)t}\, {\rm d}t. \qed
$$

\end{pf*}


\section{Conclusion}

We have extended the c.d.f.'s and moments of lattice polynomial functions to the weighted case. At first glance, this extension may appear as a
simple exercise. However, it led us to nontrivial formulas, which can be directly applied to qualitative aggregation functions such as the
Sugeno integrals and their particular cases: the weighted minima, the weighted maxima, and their ordered versions.

\section*{Appendix: Alternative proof of Eq.\ (\ref{eq:egyp})}

In this appendix we present a proof of Eq.\ (\ref{eq:egyp}) without using the c.d.f.\ of $Y_p$. The main idea of this proof is based on the fact
that the expected value $\E[g(Y_p)]$ can be expressed as an $n$-dimensional integral, namely
\begin{equation}\label{eq:NDimInt}
\E[g(Y_p)]=\int_{-\infty}^{\infty}\cdots\int_{-\infty}^{\infty}g[p(\mathbf{x})]\dd F_1(x_1)\cdots{\rm d}F_n(x_n).
\end{equation}

At first glance, the evaluation of this expression requires a difficult or intractable integration. However, as we will now see, any weighted
lattice polynomial function fulfills a remarkable decomposition formula, which will enable us to calculate $\E[g(Y_p)]$ in a straightforward
manner.

Given a weighted lattice polynomial function $p:L^n\to L$ and an index $k\in [n]$, define the weighted lattice polynomial functions
$p_k^{a}:L^n\to L$ and $p_k^{b}:L^n\to L$ as
\begin{eqnarray*}
p_k^{a}(\mathbf{x}) &:=& p(x_1,\ldots,x_{k-1},a,x_{k+1},\ldots,x_n),\\
p_k^{b}(\mathbf{x}) &:=& p(x_1,\ldots,x_{k-1},b,x_{k+1},\ldots,x_n).
\end{eqnarray*}
Then, it can be shown \cite{Marc} that
\begin{equation}\label{eq:med k}
p(\mathbf{x})=\med\big[p_k^{a}(\mathbf{x}),x_k,p_k^{b}(\mathbf{x})\big]\qquad (\mathbf{x}\in L^n).
\end{equation}
This decomposition formula expresses that, for any index $k$, the variable $x_k$ can be totally isolated in $p(\mathbf{x})$ by means of a median
calculated over the variable $x_k$ and the two functions $p_k^{a}$ and $p_k^{b}$, which are independent of $x_k$.

This interesting property leads to the following lemma.

\begin{lem}\label{lemma:egfx}
Let $p:L^n\to L$ be any weighted lattice polynomial function, let $k\in [n]$, and let $g:\overline{\R}\to\overline{\R}$ be any measurable
function of bounded variation and such that $g(-\infty)$ is finite and $\lim_{y\to\infty}g(y)[1-F_k(y)]=0$. Then
\begin{equation}\label{eq:egfx}
\E[g(Y_p)] = \E[g_a(Y_{p_k^{a}})] +\E[g_b(Y_{p_k^{b}})],
\end{equation}
where
\begin{eqnarray*}
g_a(x)&:=& \int_{-\infty}^xF_k(t)\dd g(t),\\
g_b(x)&:=& g(x)-\int_{-\infty}^xF_k(t)\dd g(t).
\end{eqnarray*}
\end{lem}

\begin{pf*}{Proof.}
Let $k\in [n]$ and fix $x_j$ for all $j\neq k$. Assume $u := p_k^{a}(\mathbf{x})$ and $v := p_k^{b}(\mathbf{x})$ are finite. The other cases can
be dealt with similarly. Then, we have $u\leqslant v$ and, by (\ref{eq:med k}),
$$
p(\mathbf{x})={\rm median}[u,x_k,v].
$$
Hence, for any measurable function $g:\overline{\R}\to\overline{\R}$ of bounded variation, we have
\begin{eqnarray*}
\int_{-\infty}^{\infty} g[p(\mathbf{x})]\dd F_k(x_k)%
&=& \int_{-\infty}^{u} g(u)\dd F_k(x_k) + \int_{u}^{v} g(x_k)\dd F_k(x_k) + \int_{v}^{\infty} g(v)\dd F_k(x_k)\\
&=& g(v)-\int_u^v F_k(t) \dd g(t)\\
&=& g_a(u)+g_b(v).
\end{eqnarray*}
Finally, integrating with respect to the other variables $x_j$, $j\neq k$, and then using (\ref{eq:NDimInt}), we get the result.\qed
\end{pf*}

Observe that formula (\ref{eq:egfx}), when considered for every index $k$ and every function $g$, completely determines the expected value of
$g(Y_p)$. Indeed, repeated applications of that formula will eventually lead to integration of transformed weighted lattice polynomial functions
all of whose variables are set to either $a$ or $b$.

\begin{pf*}{Proof of Eq.\ (\ref{eq:egyp}).}
For any fixed $S\subseteq [n]$, define recursively the sequence $\{g_S^k\}_{k=0}^n$ of functions $g_S^k:\overline{\R}\to\overline{\R}$ as
$g_S^0:=g$ and, for $k\geqslant 1$,
$$
g_S^{k}(x):=
\begin{cases}
\displaystyle{\int_{-\infty}^xF_k(t)\dd g_S^{k-1}(t)}, & \mbox{if $k\notin S$},\\
\displaystyle{g_S^{k-1}(x)-\int_{-\infty}^xF_k(t)\dd g_S^{k-1}(t)}, & \mbox{if $k\in S$}.
\end{cases}
$$
Repeated applications of Lemma~\ref{lemma:egfx} eventually lead to
$$
\E[g(Y_p)] = \sum_{S\subseteq [n]} g_S^n[p(\mathbf{e}_S)].
$$

Let us now show that
$$
g_S^{n}(z)=g_S^n(-\infty)+\int_{-\infty}^z \prod_{k\in [n]\setminus S}F_k(x)\,\prod_{k\in S}[1-F_k(x)] \dd g(x).
$$
For any $S\subseteq [n]$ and any $k\in [n]$, we have
$$
\dd g_S^{k}(x)=
\begin{cases}
F_k(x)\dd g_S^{k-1}(x), & \mbox{if $k\notin S$},\\
[1-F_k(x)]\dd g_S^{k-1}(x), & \mbox{if $k\in S$},
\end{cases}
$$
and hence
$$
\dd g_S^{n}(x)=\prod_{k\in [n]\setminus S}F_k(x)\,\prod_{k\in S}[1-F_k(x)] \dd g(x),
$$
which proves the result since $g_S^{n}(-\infty)=g(-\infty)$, if $S=N$, and $0$, otherwise.\qed
\end{pf*}


\end{document}